\documentclass[twoside,leqno,10pt]{amsart}
\usepackage{amsfonts}
\usepackage{amsmath}
\usepackage{amscd}
\usepackage{amssymb}
\usepackage{amsthm}
\usepackage{amsrefs}
\usepackage{latexsym}
\usepackage{bbm}
\numberwithin{equation}{section}

\begin{document}
\newtheorem*{theorem}{Theorem}
\newtheorem{lemma}{Lemma}
\newtheorem*{corollary}{Corollary}
\numberwithin{equation}{section}
\newcommand{\dif}{\mathrm{d}}
\newcommand{\intz}{\mathbb{Z}}
\newcommand{\ratq}{\mathbb{Q}}
\newcommand{\natn}{\mathbb{N}}
\newcommand{\comc}{\mathbb{C}}
\newcommand{\rear}{\mathbb{R}} 
\newcommand{\prip}{\mathbb{P}}
\newcommand{\uph}{\mathbb{H}}
\newcommand{\majorarc}{\mathfrak{M}}
\newcommand{\minorarc}{\mathfrak{m}}

\title{Primes in Quadratic Progressions on Average}
\author{Stephan Baier \and Liangyi Zhao}
\date{\today}
\maketitle

\begin{abstract}
In this paper, we establish a theorem on the distribution of primes in quadratic progressions on average.
\end{abstract}



\section{Introduction and Statements of Results}

It was due to Dirichlet that any linear polynomial represents infinitely many primes provided the coefficients are co-prime.  Though long been conjectured, analogous statements are not known for any polynomial of higher degree.  G. H. Hardy and J. E. Littlewood \cite{GHHJEL} conjectured that
\begin{equation} \label{hlconj}
\sum_{n \leq x} \Lambda(n^2+k) \sim \mathfrak{S}(k) x,
\end{equation}
where $\Lambda$ is the von Mangoldt function and $\mathfrak{S}(k)$ is a constant that depends only on $k$.  Their conjecture is in an equivalent but different form as in \eqref{hlconj}.  Moreover, they also gave conjectures regarding the representation of primes by any quadratic polynomial that may conceivably represent infinitely many primes.  \newline

In this paper, we aim to prove that \eqref{hlconj} holds for almost all square-free  $k \leq y$ if $x^2 (\log x)^{-A} \leq y \leq x^2$.  More in particular, we shall prove the following.

\begin{theorem} \label{mainresult}
Given $A,B>0$, we have, for $x^2 (\log x)^{-A} \leq y \leq x^2$,
\begin{equation} \label{theoeq}
\sum_{\substack{k \leq y\\ \mu^2(k)=1}} \left| \sum_{n \leq x} \Lambda(n^2+k) - \mathfrak{S}(k) x \right|^2 = O \left( \frac{yx^2}{(\log x)^B} \right),
\end{equation}
where 
\[ \mathfrak{S}(k) = \prod_{p>2} \left( 1 - \frac{\left(\frac{-k}{p} \right)}{p-1} \right) \]
with $\left( \frac{-k}{p} \right)$ being the Legendre symbol.
\end{theorem}

From the theorem, we have the following corollary.

\begin{corollary}
Given $A, B, C>0$ and $\mathfrak{S}(k)$ as defined in the theorem, we have, for $x^2 (\log x)^{-A} \leq y \leq x^2$, that 
\begin{equation} \label{coroeq}
 \sum_{n \leq x} \Lambda(n^2+k) = \mathfrak{S}(k) x + O \left( \frac{x}{(\log x)^B} \right)
\end{equation}
holds for all square-free $k$ not exceeding $y$ with at most 
$O \left( y (\log x)^{-C} \right)$ exceptions.
\end{corollary}

We note here that if we set 
\[ L(k) = \prod_{p >2} \left( 1 - \frac{ \left( \frac{-k}{p} \right)}{p} \right)^{-1}, \]
then
\[  \mathfrak{S}(k) L(k) = \prod_{p >2}  \left( 1 - \frac{\left(\frac{-k}{p} \right)}{p-1} \right) \left( 1 - \frac{ \left( \frac{-k}{p} \right)}{p} \right)^{-1} = \prod_{p >2} \frac{p^2-p-p \left(\frac{-k}{p} \right)}{p^2-p-(p-1)\left(\frac{-k}{p} \right)}. \]
Note that
\[ \frac{p^2-2p}{p^2-2p+1} \leq \frac{p^2-p-p \left(\frac{-k}{p} \right)}{p^2-p-(p-1)\left(\frac{-k}{p} \right)} \leq \frac{p^2}{p^2-1}. \]
Therefore, we have
\begin{equation} \label{boundsprod}
 \prod_{p>2} \frac{p^2-2p}{p^2-2p+1} \leq \mathfrak{S}(k) L(k) \leq \prod_{p>2} \frac{p^2}{p^2-1}.
\end{equation}
It can easily be shown that the infinite products in both the majorant and minorant of the above converge absolutely to limits that are independent of $k$.  Moreover, it is well-known that 
\begin{equation} \label{boundslk}
L(k)\neq 0 \; \mbox{and} \; L(k) \ll \log k,
\end{equation}
since $L(k)$ is the value of a Dirichlet $L$-function of modulus at most $2k$ at $s=1$.  Thus, the inequalities in \eqref{boundsprod} and \eqref{boundslk} imply that $\frak{S}(k)$ converges and
\[ \frak{S}(k) \gg \frac{1}{\log k} \gg \frac{1}{\log y} \gg \frac{1}{\log x}. \]
The above inequality shows that the main terms in \eqref{theoeq} and \eqref{coroeq} are indeed dominating for the $k$'s under consideration if $B>1$ and that we indeed have an ``almost all'' result. \newline

Our starting point is the identity
\begin{equation} \label{obofinterest}
\sum_{n \leq x} \Lambda(n^2+k) = \int_0^1 \sum_{m\leq z} \Lambda(m) e(\alpha m) \sum_{n \leq x} e(-\alpha(n^2+k)) \dif \alpha,
\end{equation}
where $z = x^2+y$. This identity is a consequence of the orthogonality relations for the function $e(z)$.  \newline
   
We use the circle method to study the question of interest and employ 
methods developed by H. Mikawa \cite{Mika} in studying the twin primes
conjecture on average.  Mikawa's result on twin primes is an improvement of earlier results of D. Wolke \cite{Wol} and A. F. Lavrik \cite{Lav1} and \cite{Lav2}. \newline

As usual in the circle method, we split the integration interval [0,1] into major arcs and minor arcs. We separate the beginning of the so-called singular series from the major arcs contribution which will give rise to the main term. We are left with the tail $\Phi(k)$ of the singular series and certain other error terms from the major arcs. The minor arcs contribution shall turn out to be an error term as well. Then we estimate the second moments over the square-free numbers $k\le y$ of all these error terms.  The second moment of $\Phi(k)$ is estimated using the classical large sieve, large sieve for real characters of Heath-Brown \cite{DRHB}, the P\'olya-Vinogradov inequality and a Siegel-Walfisz type estimate.  To estimate the second moments of the other error terms, we use Bessel's and Cauchy's inequalities together with two important Lemmas, Lemma~\ref{gallagherlemma} due to Gallagher \cite{Galla} and Lemma~\ref{mikawalemma} due to Wolke \cite{Wol} and Mikawa \cite{Mika}.  For the estimation of the minor arcs contribution we also need a standard bound for quadratic exponential sums due to Weyl. \newline

We use the following standard notations and conventions in number theory throughout paper: \newline

\noindent The symbol $p$ is reserved for primes.\newline
$e(z) = \exp (2 \pi i z) = e^{2 \pi i z}$.\newline 
$f = O(g)$ means $|f| \leq cg$ for some unspecified positive constant $c$. \newline
$f \ll g$ means $f=O(g)$. \newline
Following the general convention, we use $\varepsilon$ to denote a small positive constant which may not be the same at each occurrence.

\section{Preliminary Lemmas}

In this section, we quote lemmas that we shall need in the proofs of our theorem.  We begin with the following. \newline

\begin{lemma} [Gallagher] \label{gallagherlemma}
Let $2 < \Delta < N/2$ and $N<N'\le 2N$.  
For arbitrary complex numbers $a_n$, we have
\[ \int_{|\beta| \leq \Delta^{-1}} \left| \sum_{N<n\leq N'} a_n e(\beta n) \right|^2 \dif \beta \ll \Delta^{-2} \int_{N-\Delta/2}^{N'} \left| \sum_{\max\{ t,N \} < n \leq \min \{ t+\Delta/2, N' \}} a_n \right|^2 \dif t, \]
where the implied constant is absolute.
\end{lemma}

\begin{proof} This is Lemma 1 in \cite{Galla} in a slightly modified form.
\end{proof}

We shall also need the following lemma in our estimates of the error terms.

\begin{lemma} [Wolke, Mikawa] \label{mikawalemma}
Let
\[ \mathcal{J} (q, \Delta) = \sum_{\chi \bmod{q}} \int_N^{2N} \left| \sum_{t < n \leq t+q\Delta}^{\#} \chi(n) \Lambda(n) \right|^2 \dif t, \]
where the $\#$ over the summation symbol henceforth means that if $\chi = \chi_0$ then $\chi(n) \Lambda(n)$ is replaced by $\Lambda(n)-1$.  Let $\varepsilon$, $A$ and $B>0$ be given.  If $q \leq (\log N)^B$ and $N^{1/5+\varepsilon} \leq \Delta \leq N^{1-\varepsilon}$, then we have
\[ \mathcal{J} (q, \Delta) \ll (q\Delta)^2 N (\log N)^{-A}, \]
where the implied constant depends only on $\varepsilon$, $A$ and $B$.
\end{lemma}
\begin{proof}
This is from \cite{Mika} and is Lemma 2 there.  It can be proved using the tools in \cite{Wol}.
\end{proof}

We shall also need the following well-known inequality.

\begin{lemma} [Bessel] \label{besselineq}
Let $\phi_1$, $\cdots$, $\phi_r$ be orthonormal members of an inner product space $V$ over the complex numbers and $\xi \in V$. Then
\begin{equation}
\sum_{r=1}^R \left| \left( \xi, \phi_r \right) \right|^2 \leq \left\| \xi \right\|^2.
\end{equation}
\end{lemma}

\begin{proof}
This is a standard result.  See for example \cite{PRH} for a proof.
\end{proof}

To estimate the contribution on the minor arcs, we need the following lemma due to Weyl.

\begin{lemma} [Weyl] \label{weyllemma}
Given $x \ge 1$ and
\[ \frac{a}{q} - \frac{1}{q^2} \leq \alpha \leq \frac{a}{q} + \frac{1}{q^2}, \]
with $\gcd(a,q)=1$, we have
\begin{equation*}
\sum_{n\leq x} e (\alpha n^2) \ll \log x \left( x q^{-1/2} +(qx)^{1/2} \right),
\end{equation*}
where the implied constant is absolute.
\end{lemma}

\begin{proof}
See exercise 2 on page 215 of \cite{Bru}.
\end{proof}

We shall also need the following well-known results in analytic number theory.

\begin{lemma} [P\'olya-Vinogradov] \label{polyavino}
For any non-principal character $\chi \pmod{q}$ we have
\[ \left| \sum_{M <n \leq M+N} \chi (n) \right| \leq 6 \sqrt{q} \log q. \]
\end{lemma}

\begin{proof}
This is quoted from \cite{HIEK} and is Theorem 12.5 there.
\end{proof}

For completeness, we also quote the classical large sieve inequality for Dirichlet characters.

\begin{lemma} [Large Sieve] \label{classls}
Let $\{ a_n \}$ be an arbitrary sequence of complex numbers and $Q$, $M$, $N$ be integers with $Q,N>0$.  Then we have
\[ \sum_{q=1}^Q \frac{q}{\varphi(q)} \sideset{}{^{\star}}\sum_{\chi \; \bmod{ \; q}} \left| \sum_{n=M+1}^{M+N} a_n \chi (n) \right|^2 \ll (Q^2+N) \sum_{n=M+1}^{M+N} |a_n|^2, \]
where $\sideset{}{^{\star}}\sum$ means that the sum runs over primitive characters modulo the specified modulus only.
\end{lemma}

\begin{proof}
See for example \cite{HD}, \cite{PXG}, \cite{HM} or \cite{HM2} for the proof.
\end{proof}

We shall need the large sieve for real characters for the estimate of certain terms in the major arcs contribution.

\begin{lemma} [Heath-Brown] \label{hblemma}
Let $M$ and $N$ be natural numbers and let $a_1, \dots, a_n$ be arbitrary complex numbers.  Then
\[ \sum_{m\leq M} \left| \sum_{n \leq N} a_n \left( \frac{n}{m} \right) \right|^2 \ll (MN)^{\varepsilon} (M+N) \sum_{n\leq N} |a_n|^2, \]
for any $\varepsilon >0$, where the sums over $m$ and $n$ run over the square-free numbers.
\end{lemma}

\begin{proof}
This is Theorem 1 in \cite{DRHB}.
\end{proof}

\section{The Major Arcs}

In this section and next, we consider the contribution of the major arcs defined by
\begin{equation} \label{majorarcdef}
 \majorarc= \bigcup_{q \leq Q_1} \bigcup_{\substack{a=1 \\ \gcd(a,q)=1}}^q J_{q,a},
\end{equation}
where
\[ J_{q,a} = \left[ \frac{a}{q} - \frac{1}{qQ}, \frac{a}{q} + \frac{1}{qQ} \right], \; Q_1 = (\log x)^c, \; Q=x^{1-\varepsilon}\]
for some $c>0$ fixed and suitable . If $x$ is sufficiently
large, then $Q>Q_1$ and so the intervals $J_{q,a}$ with $q\le Q_1$ are 
disjoint. We will assume that this is the case throughout the sequel.
\newline

For $\alpha \in \majorarc$, we write
\[ \alpha = \frac{a}{q} + \beta, \; \mbox{with } |\beta| \leq \frac{1}{qQ}. \]
Let 
\[ S_1(\alpha) = \sum_{m \leq z} \Lambda(m) e(\alpha m) \; \mbox{ and } \; S_2(\alpha) = \sum_{n \leq x} e(-\alpha n^2). \]
We treat these sums in a manner similar to those treated in \cite{Mika}. \newline

We have
\begin{equation} \label{S1initial}
S_1(\alpha) = \sum_{m \leq z} \Lambda(m) e \left( \frac{a}{q} m \right) e(\beta m) = \sum_{\substack{m \leq z \\ \gcd(m,q)=1}} \Lambda(m) e \left( \frac{a}{q} m \right) e(\beta m) + O ((\log z)^2).
\end{equation}
Note that due to the presence of $\Lambda(m)$, the contribution from the terms with $\gcd(m,q) >1$ only comes from those $m$'s that are powers of primes dividing $q$ which can be absorbed into the $O$-term above.  It is also noteworthy that the implied constant in \eqref{S1initial} is absolute. \newline

It is elementary to note that if $\gcd(am,q) = 1$, we have
\begin{equation} \label{exptogauss}
  e \left( \frac{a}{q} m \right) = \frac{1}{\varphi(q)} \sum_{\chi \bmod{q}} \chi(am) \tau(\overline{\chi}),
\end{equation}
where 
$$\tau(\chi):=\sum\limits_{n=1}^q \chi(n)e\left(\frac{n}{q}\right)$$ 
is the Gauss sum and $\varphi(q)$ is the Euler $\varphi$ function.  We thus get that the first term in \eqref{S1initial} is
\begin{eqnarray*}
& & \frac{1}{\varphi(q)} \sum_{\chi \bmod{q}} \tau(\overline{\chi}) \chi(a) \sum_{m \leq z} \chi(m) \Lambda(m) e(\beta m) \\
& = & \frac{\mu(q)}{\varphi(q)} \left( \sum_{m\leq z} e(\beta m) + \sum_{m\leq z } \left( \Lambda(m) -1 \right) e (\beta m) \right) + \frac{1}{\varphi(q)} \sum_{\substack{ \chi \bmod{q} \\ \chi \neq \chi_0}} \tau(\overline{\chi}) \chi(a) \sum_{m\leq z} \chi(m) \Lambda(m) e(\beta m) \\
& = & \frac{\mu(q)}{\varphi(q)}  \sum_{m\leq z} e(\beta m) + \frac{1}{\varphi(q)} \sum_{\chi \bmod{q}} \tau(\overline{\chi}) \chi(a) \sum_{m\leq z}^{\#} \chi(m) \Lambda(m) e(\beta m) \\
& = & T_1(\alpha) + E_1(\alpha),
\end{eqnarray*}
say, where $\mu(q)$ is the M\"obius $\mu$ function and the meaning of the $\#$ over the summation symbol is the same as that in Lemma~\ref{mikawalemma}.  We arrive at
\[ S_1(\alpha) = T_1(\alpha) + E_1(\alpha) + O((\log z)^2). \]

We treat $S_2(\alpha)$ in a similar way, using \eqref{exptogauss},
\[ S_2(\alpha) = \sum_{n\leq x} e \left( -\frac{a}{q} n^2 \right) e (-\beta n^2) = \sum_{d|q} \frac{1}{\varphi(q_1^*)} \sum_{\chi \bmod{q_1^*}} \tau(\overline{\chi}) \chi(-ad^*) \sum_{\substack{n\leq x \\ \gcd(n,q)=d}} \chi^2(n^*) e(-\beta n^2), \]
where
\[ n^* = \frac{n}{d}, \; q^* =\frac{q}{d}, \; d^*= \frac{d}{\gcd(d, q^*)}, \; \mbox{and} \; q_1^* = \frac{q^*}{\gcd(d, q^*)}. \]
Hence, we get the following.
\begin{eqnarray*}
S_2(\alpha) & = & \sum_{d|q} \frac{1}{\varphi(q_1^*)} \sum_{\substack{ \chi \bmod{q_1^*} \\ \chi^2 = \chi_0}} \tau(\overline{\chi}) \chi(-ad^*) \sum_{\substack{n\leq x \\ \gcd(n,q)=d}} e(-\beta n^2) \\
 & & \hspace*{.3in} +  \sum_{d|q} \frac{1}{\varphi(q_1^*)} \sum_{\substack{ \chi \bmod{q_1^*} \\ \chi^2 \neq \chi_0}} \tau(\overline{\chi}) \chi(-ad^*) \sum_{\substack{n\leq x \\ \gcd(n,q)=d}} \chi^2(n^*) e(-\beta n^2) \\
 & = & T_2(\alpha) + E_2(\alpha),
\end{eqnarray*}
say.  Let $G= \left( \intz / q_1^* \intz \right)^*$ and $G^2= \{ g^2 : g \in G \}$.  Then it is easy to observe (see for example page 44 of \cite{HIEK}) that
\[ \sum_{\substack{ \chi \bmod{q_1^*} \\ \chi^2 = \chi_0}} \tau(\overline{\chi}) \chi(-ad^*) = \sum_{\substack{ \chi \bmod{q_1^*} \\ \chi^2 = \chi_0}} \chi(-ad^*) \sum_{b \bmod{q_1^*}} \chi(b) e \left( \frac{b}{q_1^*} \right) = \left[ G : G^2 \right] \sum_{\substack{ -ad^*b \equiv \square \bmod{q_1^*} \\ \gcd(b,q_1^*)=1}} e \left( \frac{b}{q_1^*} \right), \]
where the notation $n \equiv \square \bmod{q_1^*}$ means that $n$ is congruent to a square modulo $q_1^*$.  Moreover, we note that
\[ \sum_{\substack{ -ad^*b \equiv \square \bmod{q_1^*} \\ \gcd(b,q_1^*)=1}} e \left( \frac{b}{q_1^*} \right) = \frac{1}{\left[ G : G^2 \right]} \sum_{\substack{l=1 \\ \gcd(l,q_1^*)=1}}^{q_1^*} e \left( \frac{-\overline{ad^*}l^2}{q_1^*} \right) = \frac{1}{\left[ G : G^2 \right]} \sum_{\substack{l=1 \\ \gcd(l,q_1^*)=1}}^{q_1^*} e \left( \frac{-ad^*l^2}{q_1^*} \right), \]
upon noting that $\overline{ad^*}l^2 \equiv ad^* (\overline{ad^*}l)^2 \pmod{q_1^*}$, where $\overline{a}$ is the multiplicative inverse of $a$ modulo $q_1^*$. \newline

Consequently, we have
\[ T_2(\alpha) = \sum_{d|q} \frac{1}{\varphi(q_1^*)} \sum_{\substack{l=1 \\ \gcd(l,q_1^*)=1}}^{q_1^*} e \left( \frac{-ad^*l^2}{q_1^*} \right) \sum_{\substack{n\leq x \\ \gcd(n,q)=d}} e(-\beta n^2). \]
Furthermore, we have
\begin{equation} \label{totalmajorarc}
\begin{split}
\int_{\majorarc} \sum_{m\leq z} \Lambda(m) e(\alpha m) &\sum_{n\leq x} e(-\alpha(n^2+k)) \dif \alpha \\
&  = \int_{\majorarc} (T_1(\alpha) + E_1(\alpha) + O((\log x)^2))(T_2(\alpha)+E_2(\alpha)) e(-\alpha k) \dif \alpha.
\end{split}
\end{equation}

\section{The Singular Series}

We first consider the main term which will be given by the following
\begin{equation} \label{mainterm1}
\int_{\majorarc} T_1(\alpha) T_2(\alpha) e(-k \alpha) \dif \alpha = \sum_{q \leq Q_1} \frac{\mu(q)}{\varphi(q)} \sum_{\substack{a=1 \\ \gcd(a,q)=1}}^q e\left( - \frac{a}{q}k \right) \sum_{d|q} \frac{\mathcal{G}(a,q_1^*)}{\varphi(q_1^*)} \int_{|\beta|\leq \frac{1}{qQ}} \Pi_{q,d}(\beta) \dif \beta,
\end{equation}
where 
\[ \Pi_{q,d}(\beta)=\sum_{m \leq z} e(\beta m) \sum_{\substack{n\leq x \\ \gcd(n,q)=d}} e(-\beta n^2) e(-k\beta) \]
and
\[ \mathcal{G} (a,q_1^*) = \sum_{\substack{l=1 \\ \gcd(l,q_1^*)=1}}^{q_1^*} e \left( \frac{-ad^*l^2}{q_1^*} \right). \]

The integral on the right-hand side of \eqref{mainterm1} is well-approximated by
\begin{equation} \label{approxint}
\int_0^1 \sum_{m \leq z} e(\beta m) \sum_{\substack{n\leq x \\ \gcd(n,q)=d}} e(-\beta n^2) e(-k\beta) \dif \beta + O \left( \int_{1/(qQ)}^{1/2} \frac{1}{\beta} \left| \sum_{\substack{n\leq x \\ \gcd(n,q)=d}} e(-\beta n^2) e(-k\beta) \right| \dif \beta\right),
\end{equation}
where we have used the bound for the geometric sum over $m$ in the $O$-term above.  Applying Cauchy's inequality and Parseval's identity, the $O$-term in \eqref{approxint} is
\begin{equation} \label{errorint}
 \ll \left( \int_{1/(qQ)}^{1/2} \frac{1}{\beta^2} \ \dif \beta \right)^{\frac{1}{2}} \left( \int_0^1 \left|  \sum_{\substack{n\leq x \\ \gcd(n,q)=d}} e(-\beta n^2) e(-k\beta) \right|^2 \dif \beta \right)^{\frac{1}{2}} \ll \left( qQ\frac{x}{d} \right)^{\frac{1}{2}}.
\end{equation}
The first term in \eqref{approxint} is, by orthogonality of $e(z)$,
\begin{equation} \label{maintermint}
 \mathop{\sum_{m \leq z} \sum_{n \leq x}}_{\substack{m=n^2+k \\ \gcd(n,q)=d}} 1 = \sum_{\substack{n\leq x \\ \gcd(n,q)=d}} 1 = \sum_{\substack{ n^* \leq x/d \\ \gcd(n^*, q/d)=1}} 1 = \frac{\varphi(q/d)}{q/d} \frac{x}{d} + O \left( \varphi(q/d) \right) = \frac{\varphi(q/d)}{q} x + O \left( \varphi(q/d) \right).
\end{equation}
Now combining \eqref{mainterm1}, \eqref{approxint}, \eqref{errorint} and \eqref{maintermint}, \eqref{mainterm1} becomes
\begin{equation} \label{mainterm2}
\sum_{q \leq Q_1} \frac{\mu(q)}{\varphi(q)} \sum_{\substack{a=1 \\ \gcd(a,q)=1}}^q e\left( - \frac{a}{q}k \right) \sum_{d|q} \frac{\mathcal{G}(a,q_1^*)}{\varphi(q^*)} \left( \frac{\varphi(q/d)}{q} x + O \left( \left( qQ\frac{x}{d} \right)^{\frac{1}{2}} \right) \right).
\end{equation}
Note that $\varphi(q/d)$ is much smaller than and hence negligible in comparison with the $O$-term in \eqref{mainterm2}.\newline

Now due to the presence of $\mu(q)$ in \eqref{mainterm2}, it suffices to consider only those $q$'s that are square-free.  In that case, we have 
$d^*=d$ and $q_1^*=q^*=q/d$.  Therefore, \eqref{mainterm2} becomes
\begin{equation} \label{mainterm3}
x \sum_{q \leq Q_1} \frac{\mu(q)}{\varphi(q)q} \sum_{\substack{a=1 \\ \gcd(a,q)=1}}^q e \left( - \frac{a}{q} k \right) \sum_{d|q} \sum_{\substack{l=1 \\ \gcd(l,q/d)=1}}^{q/d} e \left( - \frac{a}{q} (ld)^2 \right) + O \left( \sqrt{xQ} (\log x)^{c_1} \right),
\end{equation}
for some fixed $c_1>0$.  It can be easily observed that
\[ \sum_{d|q} \sum_{\substack{l=1 \\ \gcd(l,q/d)=1}}^{q/d} e \left( - \frac{a}{q} (ld)^2 \right) = \sum_{r=1}^q e \left( -\frac{a}{q} r^2 \right). \]
Thus the first term in \eqref{mainterm3} becomes
\begin{equation} \label{mainterm4}
 x \sum_{q \leq Q_1} \frac{\mu(q)}{\varphi(q)q} \Sigma(q),
\end{equation}
where
\[ \Sigma(q) = \sum_{r=1}^q \sum_{\substack{a=1 \\ \gcd(a,q)=1}}^q e \left( - \frac{a}{q} (k+r^2) \right). \]
For primes $p$, we have
\begin{equation*} 
 \Sigma(p) = \sum_{r=1}^p \left( \sum_{a=1}^p e \left( - \frac{a}{p} (k+r^2) \right) -1 \right) = \sum_{\substack{r=1 \\ p|(r^2+k)}}^p p -p.
\end{equation*}
From this it follows that
\begin{equation}\label{sigmaatprime}
\Sigma(p)=\left\{ \begin{array}{llll} 0 & \mbox{ if } p=2,\\ 
p \left( \frac{-k}{p} \right) & \mbox{ if } p>2, \end{array}\right.
\end{equation}
where $\left( \frac{-k}{p} \right)$ is the Legendre symbol. Here we have used 
that $\left( \frac{-k}{p} \right)+1$ is the number of solutions to the 
congruence relation
\[ x^2 +k \equiv 0 \pmod p \]
if $p$ is an odd prime.
It can also be seen that $\Sigma(q)$ is multiplicative in the following way.  Given $q_1$ and $q_2$ with $\gcd(q_1,q_2)=1$, we have
\[ \Sigma(q_1) = \sum_{r_1=1}^{q_1} \sum_{\substack{a_1=1 \\ \gcd(a,q_1)=1}}^{q_1} e \left( - \frac{a_1}{q_1} (k+(q_2r_1)^2) \right) \]
and similarly
\[ \Sigma(q_2) = \sum_{r_2=1}^{q_2} \sum_{\substack{a_2=1 \\ \gcd(a,q_2)=1}}^{q_2} e \left( - \frac{a_2}{q_2} (k+(q_1r_2)^2) \right) \]
since, by coprimality of $q_1,q_2$, if $r_1$ and $r_2$ run over all residue classes modulo $q_1$ and $q_2$ respectively, then so do $q_2r_1$ and $q_1r_2$ respectively. Hence
\[ \Sigma(q_1)\Sigma(q_2) = \sum_{r_1=1}^{q_1} \sum_{\substack{a_1=1 \\ \gcd(a,q_1)=1}}^{q_1} \sum_{r_2=1}^{q_2} \sum_{\substack{a_2=1 \\ \gcd(a,q_2)=1}}^{q_2} e \left( f(k,a_1,a_2,q_1,q_2,r_1,r_2) \right), \]
where
\[ f(k,a_1,a_2,q_1,q_2,r_1,r_2) = -k \frac{a_1q_2+a_2q_1}{q_1q_2} - \frac{a_1q_2(q_2r_1)^2+a_2q_1(q_1r_2)^2}{q_1q_2}. \]
Note that the above is
\[ \equiv -k \frac{a_1q_2+a_2q_1}{q_1q_2} - \frac{(a_1q_2+a_2q_1)(q_1r_2+q_2r_1)^2}{q_1q_2} \pmod{1}. \]
As $a_1$ and $a_2$ run over the primitive residue classes modulo $q_1$ and $q_2$ respectively, $a_1q_2+a_2q_1$ runs over the primitive residue classes modulo $q_1q_2$; and as $r_1$ and $r_2$ run over the residue classes modulo $q_1$ and $q_2$ respectively, $r_1q_2+r_2q_1$ runs over the primitive residue classes modulo $q_1q_2$.  Therefore, we have that
\[ \Sigma(q_1)\Sigma(q_2) = \Sigma(q_1q_2). \]
In other words, $\Sigma(q)$ is multiplicative.  This fact, together with \eqref{sigmaatprime}, gives that if $q$ is square-free, then
\[ \Sigma(q) = \left\{ \begin{array}{llll} 0 & \mbox{ if } 2|q,\\ 
q \left( \frac{-k}{q} \right) & \mbox{ if } 2\nmid q,\end{array} \right. \]
where $\left( \frac{-k}{q} \right)$ is now the Jacobi symbol.  From the above, we infer that \eqref{mainterm4} is
\begin{equation} \label{mainterm5}
x \sum_{\substack{q\leq Q_1\\ 2\nmid q}} \frac{\mu(q)}{\varphi(q)} \left( \frac{-k}{q} \right) = \mathfrak{S}(k) x + O(x|\Phi(k)|),
\end{equation}
where
$$
\mathfrak{S}(k)=\sum_{\substack{q=1\\ 2\nmid q}}^{\infty} 
\frac{\mu(q)}{\varphi(q)} \left( \frac{-k}{q} \right), \; \Phi(k):= 
\sum_{\substack{q>Q_1\\ 2\nmid q}} 
\frac{\mu(q)}{\varphi(q)} \left( \frac{-k}{q}
\right).
$$
It is easy to show that the so-called singular series $\mathfrak{S}(k)$ can be rewritten as an Euler product:
\begin{equation} \label{euler}
\mathfrak{S}(k) = \prod_{p>2} \left( 1 - \frac{\left( \frac{-k}{p} \right)}{p-1} \right).
\end{equation}

We now infer from \eqref{mainterm3}, \eqref{mainterm5} and \eqref{euler} that
\begin{equation} \label{mainterm6}
\int_{\majorarc} T_1(\alpha) T_2(\alpha) e(-k\alpha) \dif \alpha 
=  \mathfrak{S}(k)x + O \left( x \left| \Phi(k) \right| +\sqrt{xQ} (\log x)^{c_1} \right).
\end{equation}

\section{The Estimate of the Second Moment of $\Phi(k)$}

In this section, we estimate the second moment over the square-free numbers $k\le y$ of the tail $\Phi(k)$ of the singular series, that is, we estimate
\[ \sum_{\substack{k \leq y\\ \mu^2(k)=1}} \left| \Phi(k) \right|^2. \]
Throughout this section, all sums over $q$ are restricted to odd $q$.  The above sum is majorized by
\begin{equation} \label{breakinto3}
\ll \sum_{k \leq y} \left| \sum_{Q_1<q \leq U} \frac{\mu(q)}{\varphi(q)} \left( \frac{-k}{q} \right) \right|^2 + \sum_{\substack{k \leq y\\ \mu^2(k)=1}} \left| \sum_{U<q\leq 2^vU} \frac{\mu(q)}{\varphi(q)} \left( \frac{-k}{q} \right) \right|^2 + \sum_{k \leq y} \left| \sum_{q>2^vU} \frac{\mu(q)}{\varphi(q)} \left( \frac{-k}{q} \right) \right|^2,
\end{equation}
with real numbers $U>Q_1$ and $v\ge 1$ to be chosen later. \newline

The first term in \eqref{breakinto3} is treated as
\begin{eqnarray*}
\sum_{k \leq y} \left| \sum_{Q_1<q \leq U} \frac{\mu(q)}{\varphi(q)} \left( \frac{-k}{q} \right) \right|^2 & = & \mathop{\sum \sum}_{Q_1 < q_1, q_2 \leq U} \frac{\mu(q_1) \mu(q_2)}{\varphi(q_1) \varphi(q_2)} \sum_{k \leq y} \left( \frac{-k}{q_1} \right) \left( \frac{-k}{q_2} \right) \\
 & \leq & y \sum_{Q_1 < q \leq U} \frac{\mu^2(q)}{\varphi^2(q)} + \mathop{\sum \sum}_{\substack{Q_1 < q_1, q_2 \leq U\\ q_1\not= q_2}} \frac{\mu(q_1) \mu(q_2)}{\varphi(q_1) \varphi(q_2)} \sum_{k \leq y} \left( \frac{-k}{q_1q_2/\gcd(q_1,q_2)} \right).
\end{eqnarray*}
The first term above is
\[ \ll \frac{y}{Q_1} \ll \frac{y}{(\log x)^{c}} \] 
and the second term, by Lemma~\ref{polyavino}, is
\[ \ll  \mathop{\sum \sum}_{\substack{Q_1 < q_1, q_2 \leq U\\ q_1\not= q_2}} \frac{\sqrt{q_1q_2}}{\varphi(q_1) \varphi(q_2)} \log (q_1q_2) \ll \log U \left( \sum_{Q_1 < q \leq U} \frac{\sqrt{q}}{\varphi(q)} \right)^2 \ll U (\log U)^{c_2}, \]
with some $c_2>0$.  Therefore, we have
\begin{equation} \label{firsttermbound}
\sum_{k \leq y} \left| \sum_{Q_1<q \leq U} \frac{\mu(q)}{\varphi(q)} \left( \frac{-k}{q} \right) \right|^2 \ll \frac{y}{(\log x)^{c}} + U (\log U)^5.
\end{equation}

To estimate the second term in \eqref{breakinto3}, we use both the classical large sieve inequality, Lemma~\ref{classls}, and the large sieve for real characters of Heath-Brown \cite{DRHB}, Lemma~\ref{hblemma}.  Using Cauchy's inequality, we have
\[ \sum_{\substack{k \leq y\\ \mu^2(k)=1}} \left| \sum_{U<q\leq 2^vU} \frac{\mu(q)}{\varphi(q)} \left( \frac{-k}{q} \right) \right|^2 \ll v \sum_{r=1}^{[v+1]} T_r, \]
where 
\[ T_r:=\sum_{\substack{ k \leq y\\ \mu^2(k)=1}} \left| \sum_{2^{r-1}U < q \leq 2^rU} \frac{\mu(q)}{\varphi(q)} \left( \frac{-k}{q} \right) \right|^2. \]
Lemma~\ref{hblemma} gives that
\[ T_r \ll (Uy)^{\varepsilon} 2^{r\varepsilon} + U^{\varepsilon-1} y^{\varepsilon+1}2^{r(\varepsilon-1)}. \]
Summing the above over $r$ with $2^rU \leq y^{1+\varepsilon}$ we get that
\begin{equation} \label{2ndtermest1}
\sum\limits_{\substack{r\ge 1 \\  2^rU \leq y^{1+\varepsilon}}} T_r
\ll (Uy)^{\varepsilon} \sum_{r \leq R} 2^{r\varepsilon} +  U^{\varepsilon-1} y^{\varepsilon+1} \sum_{r \leq R} 2^{r(\varepsilon-1)} \ll R (Uy2^R)^{\varepsilon} + U^{\varepsilon-1} y^{\varepsilon+1},
\end{equation}
where
\[ R = \left[ \log_2 \frac{y^{1+\varepsilon}}{U} \right]. \]
It is now easily observed that \eqref{2ndtermest1} is
\begin{equation} \label{2ndtermest2}
\ll y^{3 \varepsilon} + U^{\varepsilon-1} y^{\varepsilon+1}.
\end{equation}
Using the classical large sieve inequality, Lemma~\ref{classls}, we have
\begin{equation} \label{afterclassls}
T_r\ll (2^rU + y^2) \frac{1}{2^rU}.
\end{equation}
We need not worry about the primitivity of characters that is required by Lemma~\ref{classls}, since it is well-known that if $k \equiv 3 \pmod{4}$ and is square-free then $\left( \frac{-k}{q} \right)$ is primitive; if $k \equiv 1 \pmod{4}$ and is square-free then $\left( \frac{k}{q} \right)$ is primitive and $\left( \frac{-k}{q} \right)$ = $\left( \frac{-1}{q} \right)\left( \frac{k}{q} \right)$; and if $k$ is square-free and even then the Jacobi symbol is of conductor $|k/2|$.  See, for example, \S 5 of \cite{HD}.  Hence each primitive character appears at most a bounded number of times in $T_r$.  We may break the summation over $k$ into three pieces according to whether $k \equiv 1, 2$ or $3 \pmod{4}$ and then apply Lemma~\ref{classls} to each of the resulting pieces.  Summing \eqref{afterclassls} over $r$ with $y^{1+\varepsilon} < 2^rU$ and $r \leq [v+1]$, we obtain that 
\begin{equation} \label{2ndtermest3}
\sum\limits_{\substack{r\le [v+1]\\ y^{1+\varepsilon} < 2^rU}} T_r\ll v+ y^{1-\varepsilon}.
\end{equation}
Let $W =2^vU$.  For the third term in \eqref{breakinto3}, we get
\begin{equation} \label{3rdtermest1}
 \sum_{q >W} \frac{\mu(q)}{\varphi(q)} \left( \frac{-k}{q} \right) \ll \exp \left( -c \sqrt{\log W} \right)
\end{equation}
if $k\le y \ll (\log W)^{2/\varepsilon}$ by using a similar argument as in the proof of the classical Siegel-Walfisz theorem (see \cite{Bru}, Satz 3.3.3).  We note here that if $k$ is even, then the restriction of $q$ being odd on the sum over $q$ on the left-hand side of \eqref{3rdtermest1} can be removed with no change to the value of the sum, due to the fact that $\left( \frac{-k}{q} \right)$ is a character modulo $k$; and if $k$ is odd, then the sum on the left-hand side of \eqref{3rdtermest1} remains unaltered if the character $\left( \frac{-k}{q} \right)$ is replaced by the character modulo $2k$ induced by $\left( \frac{-k}{q} \right)$ and the sum over $q$ is extended over all $q$'s, both even and odd.  In both cases, a Siegel-Walfisz type argument yields the estimate in \eqref{3rdtermest1}. \newline

Now we set
\[ v = \log_2 \frac{\exp\left(y^{\varepsilon/2}\right)}{U}, \; \; U=\sqrt{y} \]
and assume without loss of generality that $\varepsilon\le 1/5$. 
Then we have
\[  v \ll y^{\varepsilon/2},\ \ \ \ \ W = \exp \left( y^{\varepsilon/2} \right),\ \ \ \ \ y= (\log W)^{2/\varepsilon},\]
and from \eqref{2ndtermest2}, \eqref{2ndtermest3} and \eqref{3rdtermest1}, we get that the sum of second and third terms in \eqref{breakinto3} is
\[ \ll v \left( y^{3 \varepsilon} + U^{\varepsilon-1} y^{\varepsilon+1} + v + y^{1-\varepsilon} \right) \ll y^{1-\varepsilon/2}. \]
Therefore, combining \eqref{firsttermbound} and the above, we get that \eqref{breakinto3} is
\begin{equation} \label{totalerrorfrommainterm}
 \ll \frac{y}{(\log x)^{c}}.
\end{equation}

\section{The Error Terms from the Major Arcs}

We consider the second moment over $k$ of the remaining terms in \eqref{totalmajorarc} term by term.  First, by Bessel's inequality, Lemma~\ref{besselineq}, we have
\begin{equation} \label{T1E2beforeest}
\begin{split}
\sum_{k \leq y} \left| \int_{\majorarc} T_1(\alpha) E_2(\alpha) e(-k\alpha) \dif \alpha \right|^2 & \ll \int_{\majorarc} \left| T_1(\alpha) E_2(\alpha) \right|^2 \dif \alpha \ll \sup_{\alpha \in \majorarc} \left| T_1 (\alpha) \right|^2 \int_{\majorarc} \left| E_2 (\alpha) \right|^2 \dif \alpha \\
&\ll z^2 \int_{\majorarc} \left| E_2 (\alpha) \right|^2 \dif \alpha.
\end{split}
\end{equation}
Now, to estimate the second factor of the above, we have
\begin{equation} \label{E2est}
\int_{\majorarc} \left| E_2 (\alpha) \right|^2 \dif \alpha = \sum_{q\leq Q_1} \sum_{\substack{a=1 \\ \gcd(a,q)=1}}^q \int_{|\beta| \leq \frac{1}{qQ}} \Omega_q(\beta) \dif \beta,
\end{equation}
where 
\[ \Omega_q(\beta):=\left| \sum_{d|q} \frac{1}{\phi(q_1^*)} \sum_{\substack{\chi \bmod{q_1^*} \\ \chi^2 \neq \chi_0}} \tau(\overline{\chi}) \chi(ad^*) \sum_{\substack{n\leq x \\ \gcd(n,q)=d }} \chi^2(n^*) e(-\beta n^2) \right|^2. \]
Applying Cauchy's inequality to $\Omega_q(\beta)$ after breaking the sum over $n$ into dyadic intervals of the form $N < n \leq N' \leq 2N \leq x$ and using the fact that $|\tau(\overline{\chi})| = \sqrt{q_1^*}$, we get that in order to estimate \eqref{E2est}, it suffices to estimate 
\begin{equation} \label{E2est2}
 (\log x)^{c_3} \sum_{q \leq Q_1} \sum_{d|q} \sum_{\substack{\chi \bmod{q_1^*} \\ \chi^2 \neq \chi_0}} \int_{|\beta| \leq \frac{1}{qQ}} \left| \sum_{\substack{N<n\leq N' \\ \gcd(n,q)=d}} \chi^2(n^*) e(-\beta n^2) \right|^2 \dif \beta,
\end{equation}
for some constant $c_3>0$. We now apply Gallagher's lemma, Lemma~\ref{gallagherlemma}, to the integral above.  We get that, if $x$ is sufficiently large, \eqref{E2est2} is majorized by
\[ \ll \frac{(\log x)^{c_3}}{Q^2} \sum_{q \leq Q_1} \frac{1}{q^2} 
\sum_{d|q} \sum_{\substack{\chi \bmod{q_1^*} \\ \chi^2 \neq \chi_0}} 
\int_{-x}^{2x} \left| \sum_{\substack{\max \{ t, N \} < n \leq \min \{ t+qQ/2, N' \} \\ \gcd(n,q)=d}} \chi^2 (n^*) \right|^2 \dif t. \]
Here we have used that $q\le Q_1=(\log x)^c$ and $Q=x^{1-\varepsilon}$.
The character sum above can be rewritten in the form
\begin{eqnarray*}
\sum\limits_{\substack{\max \{ t, N \} < n \leq \min \{ t+qQ/2, N' \} \\ \gcd(n,q)=d}} \chi^2 (n^*) &=& \sum\limits_{\substack{
\max \{ t, N \}/d < n^* \leq \min \{ t+qQ/2, N' \}/d \\ 
\gcd(n^*,q^*)=1}} \chi^2 (n^*)\\ &=& \sum\limits_{
\max \{ t, N \}/d < n^* \leq \min \{ t+qQ/2, N' \}/d} \chi' (n^*),
\end{eqnarray*}
where $\chi'$ is the (non-trivial) character modulo $q^*$ induced by 
the character $\chi^2$ mod $q_1^*$.  We now apply P\'olya-Vinogradov's estimate, Lemma~\ref{polyavino}, to the above character sum. Then, collecting all contributions, we arrive at the estimate
\begin{equation} \label{E2finalest}
\int_{\majorarc} \left| E_2 (\alpha) \right|^2 \dif \alpha \ll x Q^{-2} (\log x)^{c_4},
\end{equation}
for some fixed $c_4>0$. 
Consequently, we infer from \eqref{T1E2beforeest} and \eqref{E2finalest} that
\begin{equation} \label{T1E2est}
\sum_{k \leq y} \left| \int_{\majorarc} T_1(\alpha) E_2 (\alpha) e(-k \alpha) \dif \alpha \right|^2 \ll \frac{x^5 (\log x)^{c_4}}{Q^2}.
\end{equation}\newline

Now once again, by Bessel's inequality, Lemma~\ref{besselineq}, we have
\begin{equation} \label{E1T2beforeest}
\begin{split}
\sum_{k \leq y} \left| \int_{\majorarc} E_1(\alpha) T_2(\alpha) e(-k\alpha) \dif \alpha \right|^2 & \ll \int_{\majorarc} \left| E_1(\alpha) T_2(\alpha) \right|^2 \dif \alpha \ll \sup_{\alpha \in \majorarc} \left| T_2 (\alpha) \right|^2 \int_{\majorarc} \left| E_1 (\alpha) \right|^2 \dif \alpha \\
&\ll x^2 \int_{\majorarc} \left| E_1 (\alpha) \right|^2 \dif \alpha.
\end{split}
\end{equation}

We now need to estimate the integral in \eqref{E1T2beforeest} involving $E_1(\alpha)$.  Breaking the sum over $m$ in $E_1(\alpha)$ into dyadic intervals of the form $M<m\leq M' \leq 2M \leq z$ and applying Cauchy's 
inequality, we get that in order to estimate the integral in question, it 
suffices to estimate   
$$ 
\sum_{q\leq Q_1}  \sum_{\substack{a=1 \\ \gcd(a,q)=1}}^q \int_{|\beta| \leq \frac{1}{qQ}} \left| \frac{1}{\varphi(q)} \sum_{\chi \bmod{q}} \tau(\overline{\chi}) \chi (a) \sum_{M<m \leq M'}^{\#} \chi(m) \Lambda(m) e(\beta m) \right|^2 \dif \beta 
$$
which is precisely the same as the term $B^2$ on page 22 of \cite{Mika}.  Using the estimate for $B^2$ in \cite{Mika} which is obtained by using Lemma~\ref{gallagherlemma} and Lemma~\ref{mikawalemma}, we get that this expression is bounded by    
\begin{equation} \label{E1est}
\ll \sum_{q \leq Q_1} \frac{q}{\varphi(q)} (qQ)^{-2} \mathcal{J}(q, Q/2) + Q_1^3 Q (\log x)^2 \ll z (\log z)^{c-A},
\end{equation}
where $\mathcal{J}(q,\Delta)$ is defined in Lemma~\ref{mikawalemma}.  \newline

Now from the estimates in \eqref{E1T2beforeest} and \eqref{E1est}, we conclude that
\begin{equation} \label{E1T2finalest}
\sum_{k \leq y} \left| \int_{\majorarc} E_1(\alpha) T_2(\alpha) e(-k\alpha) \dif \alpha \right|^2 \ll x^2z (\log z)^{c-A} \ll \frac{x^4}{(\log x)^{c_5}},
\end{equation}
for any $c_5>0$. \newline

We now observe that by Cauchy's inequality and the estimates in \eqref{E2finalest} and \eqref{E1est}
\begin{equation} \label{finalerror}
\sum_{k \leq y} \left| \int_{\majorarc} E_1(\alpha) E_2(\alpha) e(-\alpha k) \dif \alpha \right|^2 \ll y \int_{\majorarc} \left| E_1(\alpha) \right|^2 \dif \alpha \int_{\majorarc} \left| E_2(\alpha) \right|^2 \dif \alpha \ll \frac{yx^3}{Q^2}.
\end{equation}

Thus, combining \eqref{T1E2est}, \eqref{E1T2finalest} and \eqref{finalerror}, we obtain that
\begin{equation} \label{majorarctotalerror}
\sum_{k \leq y} \left| \int_{\majorarc} \left( T_1(\alpha) E_2(\alpha) + T_2(\alpha) E_1(\alpha)+ E_1(\alpha) E_2(\alpha) \right) e(-k\alpha) \dif \alpha \right|^2 \ll \frac{x^5(\log x)^{c_4}}{Q^2} + \frac{x^4}{(\log x)^{c_5}} + \frac{yx^3}{Q^2}.
\end{equation}

\section{The Minor Arcs}
It still remains to consider the contribution from the minor arcs
\[ \minorarc = \left[ \frac{1}{Q}, 1+\frac{1}{Q} \right] - \majorarc, \]
where $\majorarc$ is defined in \eqref{majorarcdef}. \newline

We aim to have an estimate of the second moment over $k\le y$ of the minor arcs
contributions.  In particular, we need to estimate
\begin{equation} \label{totalminorarc}
\sum_{k\leq y} \left| \int_{\minorarc} \sum_{m \leq z} \Lambda(m) e(\alpha m) \sum_{n \leq x} e \left( - \alpha ( n^2+k) \right) \dif \alpha \right|^2.
\end{equation}
Using Bessel's inequality, Lemma~\ref{besselineq}, and 
Parseval's inequality, the above is
\begin{equation} \label{afterbessel}
\ll \int_{\minorarc} \left| S_1(\alpha) S_2(\alpha) \right|^2 \dif \alpha \ll \sup_{\alpha \in \minorarc} \left| S_2(\alpha) \right|^2 \int_0^1 \left| S_1(\alpha) \right|^2 \dif \alpha \ll \sup_{\alpha \in \minorarc} \left| S_2(\alpha) \right|^2 z \log z,
\end{equation}
where the $S_1(\alpha)$ and $S_2(\alpha)$ are the sums over $m$ and $n$ in the integrand of \eqref{totalminorarc}, respectively. \newline

We note that by Dirichlet approximation and the fact that $\alpha \in \minorarc$, 
\[ \frac{a}{q}-\frac{1}{qQ} \leq \alpha \leq \frac{a}{q} + \frac{1}{qQ} \]
for some $a$ and $q$ with $Q_1 < q \leq Q$ and $\gcd (a,q)=1$.  Therefore, we get, by Lemma~\ref{weyllemma},
\[ S_2(\alpha) \ll \log x \left( x q^{-1/2} +(qx)^{1/2} \right) \]
and hence
\[ \sup_{\alpha \in \minorarc} \left| S_2(\alpha) \right|^2 \ll (\log x)^2 \left( x^2 Q_1^{-1} + Qx \right). \]
Thus we infer that \eqref{afterbessel}, and hence \eqref{totalminorarc} is majorized by
\begin{equation} \label{minorarcbound}
\ll \frac{x^4}{(\log x)^{c-3}} + (\log x)^3 Qx^3.
\end{equation}

\section{Proof and Discussion of the Theorem}
Using Cauchy's inequality, 
combining \eqref{totalmajorarc}, \eqref{mainterm6}, \eqref{totalerrorfrommainterm}, \eqref{majorarctotalerror} and \eqref{minorarcbound} and recalling $Q=x^{1-\varepsilon}$, 
we obtain the theorem. \newline

We would like to note that the result could be improved substantially if we assume the Generalized Riemann Hypothesis (GRH) for Dirichlet $L$-functions. Indeed, under GRH, Lemma 2 holds for a much larger range of $q$, namely for $q \leq x^{\delta_1}$ with some positive $\delta_1$. Moreover, a better majorant for $\mathcal{J}(q,\Delta)$ would be true under the GRH.  More precisely, we would have in Lemma~\ref{mikawalemma} a saving of a positive power of $N$ rather than a saving of a power of logarithm as we have now.  These improvements would imply that the theorem holds in a much wider $y$-range, namely $x^{2-\delta_2}\le y\le x^2$, for some positive $\delta_2$. \newline
  
{\bf Acknowledgements.}
This paper was written when the first and second-named authors held postdoctoral fellowships at the Department of Mathematics and Statistics at Queen's University and the Department of Mathematics at the University of Toronto, respectively.   The authors wish to thank these institutions for their financial support.  Moreover, this work was started when the second-named author was visiting {\it La Centre de Recherches Math\'ematiques} (CRM) in {\it Universit\'e de Montr\'eal} as a guest researcher during the Theme Year 2005-2006 in Analysis in Number Theory.  He would like to thank the CRM for their financial support and warm hospitality during his pleasant stay in Montreal.

\bibliography{biblio}
\bibliographystyle{amsxport}

\vspace*{.7cm}

\noindent Department of Mathematics and Statistics, Queen's University \newline
University Ave, Kingston, ON K7L 3N6 Canada \newline
Email: {\tt sbaier@mast.queensu.ca} \newline

\noindent Department of Mathematics, University of Toronto \newline
40 Saint George Street, Toronto, ON M5S 2E4 Canada \newline
Email: {\tt lzhao@math.toronto.edu}
\end{document}